\documentclass[reqno]{amsart}

\usepackage{euscript}
\usepackage{hyperref}
\usepackage{graphicx}

\begin{document}

\title[Two representations of the fundamental
group]{Two representations of the fundamental group and invariants
of lens spaces}

\author{E.~V.~Martyushev}

\address{Southern Ural State University, 76 Lenin avenue,
454080 Chelyabinsk, Russia}

\email{mev@susu.ac.ru}

\thanks{This work has been partially supported by
Russian Foundation for Basic Research, Grant no. 01-01-00059}

\begin{abstract}

This article is a continuation of work on construction and
calculation various of modifications of invariant based on the use
Euclidean metric values attributed to elements of manifold
triangulation. We again address the well investigated lens spaces
as a standard tool for checking the nontriviality of topological
invariants.

\end{abstract}

\maketitle

\section{Introduction}\label{sec:intro}

Let us shortly remind the main constructions of
works~\cite{1},~\cite{2},~\cite{4}.

Let there be given the universal covering of an oriented
three-dimensional manifold, considered as a simplicial complex, in
which every simplex is mapped into the three-dimensional Euclidean
space (so that the intersections of their images in $\mathbb{R}^3$
are possible). To each edge is assigned thus a Euclidean length,
and each tetrahedron has a sign ``plus'' or ``minus'' depending on
orientation of its image in $\mathbb {R}^3$; the defect angle
(minus algebraic sum  of dihedral angles) at each edge is equal to
zero modulo $2\pi$. Then we consider infinitesimal translations of
vertices of the complex which lead to infinitesimal changes of
edge lengths, while these latter lead, in their turn, to the
infinitesimal changes of defect angles. Thus, the following
sequence of vector spaces arises:
\begin{equation}\label{eq:acycl1}
0\xrightarrow{{}}\mathfrak{e}_3
\xrightarrow{{C}}(dx)\xrightarrow{{B}}(dl)\xrightarrow{{A}}(d\omega
)\xrightarrow{{B^T}}( \ldots )\xrightarrow{{C^T}}( \ldots
)\xrightarrow{{}}0.
\end{equation}
Here $\mathfrak{e}_3$ is the Lie algebra of motions of
three-dimensional Euclidean space; $(dx)$ is the space of column
vectors of vertex coordinates differentials; $(dl)$ is the space
of column vectors of edge length differentials; $(d\omega)$ is the
space of column vectors of differentials of defect angles. The
superscript $T$ means matrix transposition. Matrices $A$, $B$ and
$C$ consist of partial derivatives. One can see~\cite{4} that the
sequence~(\ref{eq:acycl1}) is an algebraic complex (superposition
of two next matrices is equal to zero). Besides, it is possible to
show~\cite{4} that for manifolds with finite fundamental group
this complex is acyclic.

For lens spaces we consider two sets  of representations of
fundamental group (recall that $\pi_1(L(p,q)) = \mathbb{Z}_p$):
\begin{subequations}
\begin{align}
 f_k \colon \mathbb{Z}_p \to & \{ \text{group of motions of
 3-dimensional Euclidean space} \},\label{eq:fk} \\
 g_j \colon \mathbb{Z}_p \to & \{ \text{group of automorphisms
 of vector spaces entering in~(\ref{eq:acycl1})} \},\label{eq:gj} \\
 & j,k = 0, \ldots, p - 1. \notag
\end{align}
\end{subequations}
The realization of representations $f_k$ was described in detail
in~\cite{2} (there they were denoted $\varphi_k$). As to
representations $g_j$, we will describe them below when we
consider the structure of spaces entering in
complex~(\ref{eq:acycl1}).

\section{Calculation of invariants for $L(p,q)$}\label{sec:calc}

\begin{figure}
\centering
\includegraphics[width=5.5cm,height=8cm]{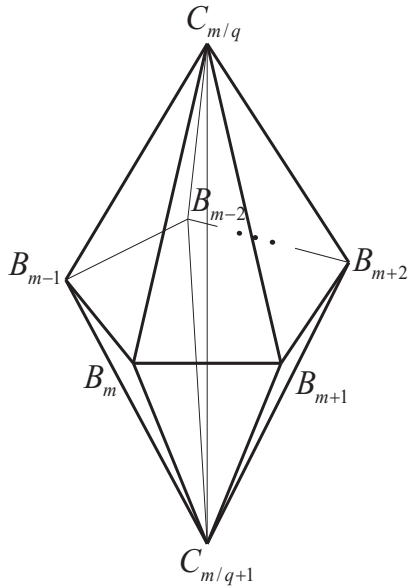}
\caption{Simplicial complex for lens space (all indexes are taken
modulo $p$)} \label{fig:complex}
\end{figure}

In figure~\ref{fig:complex} we show the simplicial complex
(bipyramid) corresponding to the lens space $L(p,q)$. Let index
$m$ (in figure) vary from $0$ to $p - 1$. Then gluing together $p$
bipyramids numbered by $m$ by identical vertices, edges and
$2$-simplices we obtain an image of simplicial complex for
$\widetilde{L(p,q)} = S^3$ in $\mathbb{R}^3$.

Every element of group $\mathbb{Z}_p$ can be represented in the
group of motions of three-dimensional Euclidean space as a
rotation through angle $2\pi k / p$ around some axis which we
denote $x_0$. The integer parameter $k$ numbers the characters of
these representations. Thus, in the right Cartesian coordinate
system $O x_1 x_{p - 1} x_0$ the coordinates of vertices of
simplicial complex look like (cf.~\cite{2}):
\begin{align}\label{eq:coords}
 B_m & \left( \rho \cos \dfrac{{2\pi mk}} {p}, \rho \sin
 \dfrac{{2\pi mk}} {p}, 0  \right), \\
 C_m & \left( \sigma \cos
 \left(\alpha + \dfrac{{2\pi q mk}} {p}\right), \sigma \sin
 \left(\alpha + \dfrac{{2\pi q mk}} {p}\right), s \right),
\end{align}
where $\alpha$, $\rho$, $\sigma$, $s$ are some parameters that
fix, to within the motions of $\mathbb{R}^3$, the location of
vertices and hence all geometrical values of the complex such as
edge lengths, volumes of tetrahedra, dihedral angles, etc. For
instance, a simple calculation shows that the oriented volumes of
tetrahedra are determined in this case by the formula:
\begin{equation}\label{eq:volumes}
V_{C_0 C_1 B_{m + 1} B_m} = \dfrac{4} {6} \rho \sigma s \sin
\dfrac{\pi k}{p} \sin \dfrac{\pi q k}{p} \sin \left( {\alpha +
\dfrac{{\pi (q - 1 - 2m)k}} {p}} \right).
\end{equation}

Let us write out now more explicitly some of vector spaces
entering in~(\ref{eq:acycl1}) (due to the symmetry of the complex,
it is enough for us to deal only with the first three of them):
\begin{equation}\label{eq:spaces}
\mathfrak{e}_3 =
 \begin{pmatrix}
   {d\varphi _1 }  \\
   {d\varphi _{p - 1} }  \\
   {d\varphi _0 }  \\
   {dx_1}  \\
   {dx_{p - 1}}  \\
   {dx_0}
 \end{pmatrix},
\quad (dx) =
 \begin{pmatrix}
   {(dx)_0 }  \\
   {(dx)_1 }  \\
    \vdots   \\
   {(dx)_{p - 1} }
 \end{pmatrix},
\quad (dl) =
 \begin{pmatrix}
   {(dl)_0 }  \\
   {(dl)_1 }  \\
    \vdots   \\
   {(dl)_{p - 1} }
 \end{pmatrix},
\end{equation}
where $dx_{\ldots}$ and $d\varphi_{\ldots}$ are infinitesimal
translations and rotations in the described above coordinate
system $O x_1 x_{p - 1} x_0$.

Every element of group $\mathbb{Z}_p$ acts on spaces $(dx)$ and
$(dl)$ by cyclic shifts as follows: it sends subspaces $(dx)_m$
and $(dl)_m$ to $(dx)_{jm}$ and $(dl)_{jm}$ respectively. Index
$j$ changes from $0$ to $p - 1$. So, these actions completely
determine representations $g_j$.

Subspaces $(dx)_m$ and $(dl)_m$ look like (all indexes are taken
modulo $p$):
\begin{equation}\label{eq:subspaces}
(dx)_m  =
 \begin{pmatrix}
   {dx_{B_m } }  \\
   {dy_{B_m } }  \\
   {dz_{B_m } }  \\
   {dx_{C_{m/q} } }  \\
   {dy_{C_{m/q} } }  \\
   {dz_{C_{m/q} } }
 \end{pmatrix},
\quad (dl)_m  =
 \begin{pmatrix}
   {dl_{B_m B_{m + 1} } }  \\
   {dl_{B_m C_{m/q} } }  \\
   {dl_{B_{m + 1} C_{m/q} } }  \\
    \vdots   \\
   {dl_{B_{m - 1} C_{m/q} } }  \\
   {dl_{C_{m/q} C_{m/q + 1} } }
 \end{pmatrix}.
\end{equation}

In general, matrices $A$, $B$ and $C$ specifying the mappings of
vector spaces in complex~(\ref{eq:acycl1}) are not block-diagonal.
In order to make them such we perform unitary change-of-basis
transformations in all spaces entering in the complex. Our goal is
to split~(\ref{eq:acycl1}) into the direct sum of subcomplexes
each of which would correspond to the character of some
representation $g_j$ (see~(\ref{eq:gj})). The unitarity of
transformations is necessary to preserve the symmetric
form~(\ref{eq:acycl1}) of the complex. Let $\varepsilon$ denote
the primitive root of unity of degree $p$. Then
\begin{itemize}
 \item the change-of-basis matrix in space $\mathfrak{e}_3$ is
\begin{equation}\label{eq:matrixU1}
U_1 =
 \begin{pmatrix}
   {\dfrac{{\sqrt 2 }}
{2}} & {\dfrac{{\sqrt 2 }} {2}i} & 0 & {} & {} & {}  \\
   {\dfrac{{\sqrt 2 }}
{2}i} & {\dfrac{{\sqrt 2 }} {2}} & 0 & {} & {} & {}  \\
   0 & 0 & 1 & {} & {} & {}  \\
   {} & {} & {} & {\dfrac{{\sqrt 2 }}
{2}} & {\dfrac{{\sqrt 2 }} {2}i} & 0  \\
   {} & {} & {} & {\dfrac{{\sqrt 2 }}
{2}i} & {\dfrac{{\sqrt 2 }} {2}} & 0  \\
   {} & {} & {} & 0 & 0 & 1
 \end{pmatrix};
\end{equation}
\item the change-of-basis matrix in space $(dx)$ is
\begin{equation}\label{eq:matrixU2}
U_2 = \dfrac{1} {{\sqrt p }}
 \begin{pmatrix}
   {1H_0 } & {1H_0 } &  \cdots  & {1H_0 }  \\
   {1H_1 } & {\varepsilon ^k H_1 } &  \cdots  & {\varepsilon ^{(p - 1)k} H_1 }  \\
    \vdots  &  \vdots  & {} &  \vdots   \\
   {1H_{p - 1} } & {\varepsilon ^{(p - 1)k} H_{p - 1} } &  \cdots  &
   {\varepsilon ^{(p - 1)^2 k} H_{p - 1} }
 \end{pmatrix},
\end{equation}
where
\begin{equation}\label{eq:matrixHm}
H_m = U_1
 \begin{pmatrix}
   {\varepsilon ^{ - mk} } & 0 & 0 & {} & {} & {}  \\
   0 & {\varepsilon ^{mk} } & 0 & {} & {} & {}  \\
   0 & 0 & 1 & {} & {} & {}  \\
   {} & {} & {} & {\varepsilon ^{ - mk} } & 0 & 0  \\
   {} & {} & {} & 0 & {\varepsilon ^{mk} } & 0  \\
   {} & {} & {} & 0 & 0 & 1
 \end{pmatrix};
\end{equation}
\item the change-of-basis matrix in space $(dl)$ is
\begin{equation}\label{eq:matrixU3}
U_3 = \dfrac{1} {{\sqrt p }}
 \begin{pmatrix}
   {1I } & {1I } &  \cdots  & {1I }  \\
   {1I } & {\varepsilon ^k I } &  \cdots  & {\varepsilon ^{(p - 1)k} I }  \\
    \vdots  &  \vdots  & {} &  \vdots   \\
   {1I } & {\varepsilon ^{(p - 1)k} I } &  \cdots  &
   {\varepsilon ^{(p - 1)^2 k} I }
 \end{pmatrix},
\end{equation}
where $I$ is the identical matrix of size $(p + 2)\times (p + 2)$.
\end{itemize}

Now, in an obvious way,
\begin{align*}
C \to U_2^{\dag} C U_1,
\\
B \to U_3^{\dag} B U_2,
\\
A \to U_3^{\dag} A U_3,
\end{align*}
and each of matrices $A$, $B$ and $C$ becomes block-diagonal (the
rearrangement of rows or/and columns can be required). Thus, the
acyclic complex~(\ref{eq:acycl1}) splits into the following direct
sum of $p$ acyclic subcomplexes:
\begin{subequations}\label{eq:acycl2}
\begin{align}
0\xrightarrow{{}}
 \begin{pmatrix}
   {d\varphi _j } \\
   {dx_j}
 \end{pmatrix}
\xrightarrow{\sqrt{p}\,{C_j }}(dx)_j\xrightarrow{{B_j }} &
(dl)_j\xrightarrow{{A_j }}(d\omega )_j\xrightarrow{{B_j^{\dag} }}(
\ldots )_j\xrightarrow{{\sqrt{p}\,C_j^{\dag} }}( \ldots
)_j\xrightarrow{{}}0,
\\
& j = 0, \pm 1 \bmod p, \notag
\\
0 \xrightarrow{{}}(dx)_j\xrightarrow{{B_j }} &
(dl)_j\xrightarrow{{A_j }}(d\omega )_j\xrightarrow{{B_j^{\dag} }}(
\ldots )_j\xrightarrow{{}}0,
\\
& j = 2, \ldots , p - 2. \notag
\end{align}
\end{subequations}
The superscript $\dag$ means operation of Hermitian conjugation.

We define the torsions of complexes~(\ref{eq:acycl2}) by formulas:
\begin{equation}\label{eq:torsions}
\EuScript{T}_j =
  \begin{cases}
 p^{- 2}\bigl|(\det C_j {\vert}_{{\overline {\mathcal D}}_j} )\bigr| ^{ - 2} \,
 \bigl|(\det {}_{{\mathcal D}_j} {\vert} B_j {\vert}_{{\overline
 {\mathcal C}}_j} )\bigr| ^2 \,(\det {}_{{\mathcal C}_j} {\vert} A_j )^{ -
 1}, & j = 0, \pm 1 \bmod p, \\[.4 cm]
 \bigl|(\det {}_{{\mathcal D}_j}{\vert} B_j
 {\vert}_{{\overline {\mathcal C}}_j} )\bigr| ^2 \,(\det
 {}_{{\mathcal C}_j}{\vert} A_j )^{ - 1}, & j = 2, \ldots,
 p - 2.
  \end{cases}
\end{equation}
Here ${\mathcal C}_j$ is a maximal subset of edges such that the
restriction of matrix $A_j$ onto its corresponding subspace of
$(dl)$ (we denote it as ${}_{{\mathcal C}_j}{\vert} A_j$) is
nondegenerate. ${\overline {\mathcal C}}_j$ is the complement of
${\mathcal C}_j$. Sets ${\mathcal D}_j$ and ${\overline {\mathcal
D}}_j$ consist of vertices and are defined similarly.

On the base of torsions~(\ref{eq:torsions}), one can
construct~\cite{4} invariant values not changing under the Pachner
moves $2 \leftrightarrow 3$ and $1 \leftrightarrow 4$. In our
case, these invariants take the form:
\begin{equation}\label{eq:invs1}
 I_{j,1}(L(p,q))  = \left\{ \EuScript{T}_j \dfrac{{l_{B_0 B_1}^2
 l_{C_0 C_1}^2 \prod\limits_{m = 0}^{p - 1}{l_{B_m C_0}^2 } }}
 {{\prod\limits_{m = 0}^{p - 1}{6V_{C_0 C_1 B_{m + 1} B_m}} }} \colon k = 1, \ldots,
 \left\lfloor {\dfrac{p} {2}} \right\rfloor \right\}, \qquad j = 0, \ldots,
 p - 1.
\end{equation}
Here $l_{\ldots}$ are usual Euclidean edge lengths, $V_{\ldots}$
are volumes of tetrahedra given by~(\ref{eq:volumes}). We use
definition $I_{j, 1}$ to emphasize that, in contrast with
work~\cite{3}, here we deal with nontrivial representations $f_k$.
The calculations implied by formula~(\ref{eq:invs1}) are as
follows. Using acyclicity of complexes, we choose subsets
${\mathcal C}_j$ and ${\mathcal D}_j$. Then we calculate the
corresponding minors of matrices $A_j$, $B_j$ and $C_j$. In
particular, we use the known expressions for partial derivatives
$\dfrac{\partial \omega}{\partial l}$ (see~\cite{2} for details).
After this, we find the torsions according to~(\ref{eq:torsions})
and substitute them into~(\ref{eq:invs1}).

The computer calculations which have been carried out for first
several $p$'s allow us to suppose that the general formula for
invariants of lens spaces~(\ref{eq:invs1}) looks as follows:
\begin{subequations}\label{eq:invs2}
\begin{align}
 I_{0,1}(L(p, q)) & = \left\{ \dfrac{(- 1)^{p - 1}}{p^4} \Delta_k^4 \colon k = 1, \ldots,
 \left\lfloor {\dfrac{p} {2}} \right\rfloor \right\}, \\
 I_{\pm 1,1}(L(p, q)) & = \left\{ \dfrac{(- 1)^{p - 1}}{p^4} \Delta_k^2 \Delta_{2k}^2
 \colon k = 1, \ldots, \left\lfloor {\dfrac{p} {2}} \right\rfloor \right\}, \\
 I_{j,1}(L(p,q)) & = \left\{ (- 1)^p \Delta_{(j - 1)k}^2 \Delta_{jk}^2 \Delta_{(j + 1)k}^2
 \colon k = 1, \ldots, \left\lfloor {\dfrac{p} {2}} \right\rfloor
 \right\},
\end{align}
\end{subequations}
where we have denoted
\begin{equation}\label{eq:Deltam}
\Delta_m  = 4\sin \dfrac{{\pi m}} {p}\sin \dfrac{{\pi q m}}{p}.
\end{equation}

\section{Discussion}\label{sec:disc}

Let us note once again that the value $\Delta_m$ determined
in~(\ref{eq:Deltam}) is nothing but the module of $m$th component
of the Reidemeister torsion of $L(p,q)$. It is of course very
intriguing to calculate invariants of such kind for other 3- and
4-manifolds.

In the future we plan to investigate similar invariants based on
$SL(2)$-solution of the pentagon equation~\cite{6},~\cite{7}.
After that, we plan to move on from lens spaces to manifolds with
noncommutative or/and infinite fundamental group.

\subsection*{Acknowledgements}

I am grateful to I.~G.~Korepanov for proposing me this problem and
improvements in the initial text of this Letter.

\bibliographystyle{amsplain}

\end{document}